\documentclass[]{elsarticle}
\usepackage{latexsym}
\usepackage{amsmath}
\usepackage{multirow}
\usepackage{graphicx}
\usepackage{amsfonts}
\usepackage{url}
\usepackage[ruled,linesnumbered,noend]{algorithm2e}
\usepackage{bm} 
\newcommand{\mb}[1]{\bm{#1}}

\DeclareMathOperator*{\argmin}{arg\,min}
\newcommand{\MA}{\bm{A}}
\newcommand{\Vb}{\bm{b}}
\newcommand{\Vbeta}{\bm{g}^c}
\newcommand{\Vg}{\bm{g}}

\newcommand{\Vp}{\bm{p}}
\newcommand{\Vphi}{\bm{g}^f}
\newcommand{\Vx}{\bm{x}}

\begin{document}
\begin{frontmatter}
\title{Active Set Expansion Strategies in MPRGP Algorithm}
\author[ugn,kam]{J.~Kru\v{z}\'ik\corref{cor}}
\ead{jakub.kruzik@ugn.cas.cz}

\author[kam,ugn]{D.~Hor\'ak}
\ead{david.horak@vsb.cz}

\author[km,kam,ugn,enet]{M.~\v Cerm\'ak}
\ead{martin.cermak@vsb.cz}

\author[km]{L.~Posp\'i\v{s}il}
\ead{lukas.pospisil@vsb.cz}

\author[ugn,kam]{M.~Pecha}
\ead{marek.pecha@vsb.cz}

\cortext[cor]{Corresponding author}
\address[ugn]{Institute of Geonics of the Czech Academy of Sciences, Ostrava, Czech Republic}
\address[kam]{Department of Applied Mathematics, FEECS, V\v SB-TU Ostrava, Ostrava, Czech Republic}
\address[km]{Department of Mathematics, Faculty of Civil Engineering, V\v SB-TU Ostrava, Ostrava, Czech Republic}
\address[enet]{ENET Centre, V\v SB-TU Ostrava, Ostrava, Czech Republic}

\begin{abstract}
The paper investigates strategies for expansion of active set that can be employed by the MPRGP algorithm. The standard MPRGP expansion uses a projected line search in the free gradient direction with a fixed step length. Such a scheme is often too slow to identify the active set, requiring a large number of expansions. We propose to use adaptive step lengths based on the current gradient, which guarantees the decrease of the unconstrained cost function with different gradient-based search directions. Moreover, we also propose expanding the active set by projecting the optimal step for the unconstrained minimization. Numerical experiments demonstrate the benefits of our expansion step modifications on two benchmarks -- contact problem of linear elasticity solved by TFETI and machine learning problems of SVM type, both implemented in PERMON toolbox.
\end{abstract}

\begin{keyword}
MPRGP \sep active set \sep expansion step \sep PERMON
\end{keyword}
\end{frontmatter}

\section{Introduction}
Quadratic programming (QP) problems with bound or box constraints can be solved by the MPRGP (Modified Proportioning with Reduced Gradient Projections) algorithm developed by Dostal \cite{Dos-book-09}, which belongs among the active set based methods.
The algorithm has been proven to enjoy an R-linear rate of convergence given by the bound on the spectrum of the Hessian matrix.
In each iteration, MPRGP performs one of three types of steps - unconstrained minimization, expansion, and proportioning.
The unconstrained minimization is typically performed by a conjugate gradient (CG) step and we will assume this for the rest of the article.
The active set is expanded by the expansion step, which consists of a maximal feasible unconstrained minimization, in our case a partial CG step to the bound/box, followed by a feasible fixed step length line search.
Finally, a proportioning step designed to reduce the active set consists of a steepest descent step in the direction of a chopped gradient.

This paper deals with the modification of the expansion step. In the original version, the theory supports fixed step length expansion step taking its value from zero to the two divided by the norm of the Hessian.
However, many numerical experiments demonstrate that this fixed step length can result in a large number of expansion steps.
To reduce the number of expansion step we propose several alternatives to the original expansion.
Taking into account the situation in the current iteration, we present two adaptive step lengths.
Additionally, we provide a comparison of the various choices for the search direction.
Both, the step lengths and the search directions are based on the current gradient splitting.
All of these expansion steps do a partial CG step to the bound/box.
A natural idea how to expand the active set is to perform the full CG step with a subsequent projection onto the feasible set.
Therefore, we propose a projected CG step as another variant of the expansion.

The benefits of our new approaches are documented on two model benchmarks - TFETI (Total Finite Element Tearing and Interconnecting) applied to a contact problem of mechanics and a sequence of machine learning problems solved by SVMs (Support Vector Machines).
PERMON (Parallel, Efficient, Robust, Modular, Object-oriented, Numerical) \cite{permon-www} toolbox, was used for the numerical experiments.

The paper is divided as follows. Section \ref{sec:mprgp} describes the MPRGP algorithm. Our modifications of the expansion step are presented in Section \ref{sec:expansion}. Section \ref{sec:soft} briefly introduces the employed software and methods used in the benchmarks. The numerical experiments are presented in Section \ref{sec:num}. Finally, we draw our conlusions in Section \ref{sec:conc}. The complete results of the numerical experiments are provided in \ref{sec:appendix}.

\section{The MPRGP Algorithm}
\label{sec:mprgp}
MPRGP  \cite{Dos-book-09} represents
an efficient algorithm for the solution of convex QP with box constraints, i.e. for minimizing quadratic functional subject to constrains
\begin{equation}
  \argmin_{\Vx} f(\Vx) = \argmin_{\Vx} \frac{1}{2}\bm{x}^{T}\bm{A}\bm{x}-\bm{x}^{T}\bm{b}\quad\text{s.t.}\quad \bm{l} \le \bm{x} \le \bm{u},
  \label{eq:MPGPprob}
\end{equation}
where $f(\Vx)$ is the cost function, $\bm{A} \in \mathbb{R}^{n\times n}$ is positive semi-definite Hessian, $\bm{x}$ is the solution, $\bm{b}$ is the right hand side, $\bm{l}$ and $\bm{u}$ is the lower respectively the upper bound.

To describe the algorithm we first have to define a gradient splitting.
Let $\bm{g}=\bm{A}\bm{x}-\bm{b}$ be the gradient.
Then we can define a component-wise (for $j \in \{1,2,\dots,n\})$ gradient splitting which is computed after each gradient evaluation. The free gradient is defined as
\begin{flalign*}
  g_j^f = \begin{cases}
  	0 \quad &\text{if}\quad x_j=l_j \quad\text{or}\quad x_j=u_j,\\
    g_j \quad &\text{otherwise}.
    \end{cases}
\end{flalign*}
The reduced free gradient is
\begin{flalign*}
   g_j^r = \begin{cases} 0 \quad&\text{if}\quad x_j=l_j \quad\text{or}\quad x_j=u_j,\\
  \min\left(\frac{x_j-l_j}{\overline{\alpha}},g_j\right) \quad &\text{if}\quad l_j < x_j <  u_j \quad\text{and}\quad g_j> 0,\\
  \max\left(\frac{x_j-u_j}{\overline{\alpha}},g_j\right)\quad &\text{if}\quad l_j < x_j <  u_j \quad\text{and}\quad g_j \le 0,
  \end{cases}
\end{flalign*}
where $\overline{\alpha} \in (0,2||\bm{A}||^{-1}]$ is used as an appriory chosen fixed step length in the expansion step.
Effectively, $\Vphi$ is the gradient on the free set and $\Vg^r$ is the free gradient that is reduced such that a step in its opposite direction with the step length $\overline{\alpha}$ does not leave the feasible set $\Omega = \{\bm{x}: \bm{l} \le \bm{x} \le \bm{u} \}$. A step in either of these direction can expand the active set, but cannot reduce it.

The chopped gradient is defined as
\begin{flalign*}
  g_j^c = \begin{cases} 0 \quad&\text{if}\quad l_j < x_j < u_j,\\
  \min(g_j,0) \quad&\text{if}\quad  x_j = l_j,\\
  \max(g_j,0) \quad&\text{if}\quad x_j = u_j.
  \end{cases}
\end{flalign*}
A step in the direction opposite of $\Vbeta$ may reduce the active set, but cannot expand it.

The next ingredient is the projection onto the feasible set $\Omega$ which is defined as 
\begin{displaymath}
  \left[P_{\Omega}(\bm{x})\right]_j = \min(u_j,\max(l_j,x_j)).
\end{displaymath}

Finally, the projected gradient is defined as $\bm{g}^{P}= \bm{g}^f + \bm{g}^c$. Its norm decrease is the natural stopping criterion of the algorithm.

These are all the necessary ingredients to summarise MPRGP in Algorithm \ref{alg:mprgp}.

\begin{algorithm}[htb]
    \SetStartEndCondition{ }{}{}%
    \DontPrintSemicolon
    \SetKwProg{Fn}{def}{\string:}{}
    \SetKwFunction{Range}{range}
    \SetKw{KwTo}{in}\SetKwFor{For}{for}{\string:}{}%
    \SetKwIF{If}{ElseIf}{Else}{if}{:}{elif}{else:}{}%
    \SetKwFor{While}{while}{:}{fintq}%
    \AlgoDontDisplayBlockMarkers\SetAlgoNoEnd\SetAlgoNoLine%
    \KwIn{$\MA$, $\Vx^{0} \in \Omega$, $\Vb$, $\Gamma > 0$, $\overline{\alpha} \in (0,2||\MA||^{-1}]$}
    $\Vg = \MA \Vx_{0} - \Vb$, $\Vp = \Vphi(\Vx^0)$, $k = 0$\;
    \While{$||\Vg^P||$ is not small}{
        \If{$||\Vbeta||^2 \le \Gamma^2 ||\Vg^f||^2$}{
            $\alpha_{f} = \max \{ \alpha_{cg}: \Vx^k - \alpha_{cg} \Vp \}$\;    
            $\alpha_{cg} = \Vg^T \Vp / \Vp^T \MA \Vp$\;
            \If{$\alpha_{cg} \le \alpha_{f}$}{
              CG() - Algorithm \ref{alg:cg} \;
            }\Else{
              Expansion() - Algorithm \ref{alg:exp};
            }
        }\Else{
            Proportioning() - Algorithm \ref{alg:prop};
        }
        $k=k+1$\;
    }
    \KwOut{$\Vx^{k}$}
    
    \caption{MPRGP}
    \label{alg:mprgp}
\end{algorithm}

\begin{algorithm}[htb]
    \SetStartEndCondition{ }{}{}%
    \DontPrintSemicolon
    \SetKwProg{Fn}{def}{\string:}{}
    \SetKwFunction{Range}{range}
    \SetKw{KwTo}{in}\SetKwFor{For}{for}{\string:}{}%
    \SetKwIF{If}{ElseIf}{Else}{if}{:}{elif}{else:}{}%
    \SetKwFor{While}{while}{:}{fintq}%
    \AlgoDontDisplayBlockMarkers\SetAlgoNoEnd\SetAlgoNoLine%
    $\Vx^{k+1} = \Vx^{k} -\alpha_{cg}\Vp$\;
    $\Vg = \Vg - \alpha_{cg} \MA \Vp$\;
    $\beta = \Vp^T \MA \Vphi / \Vp^T \MA \Vp$\;
    $\Vp = \Vphi -\beta \Vp$\;
    \caption{CG}
    \label{alg:cg}
\end{algorithm}

\begin{algorithm}[htb]
    \SetStartEndCondition{ }{}{}%
    \DontPrintSemicolon
    \SetKwProg{Fn}{def}{\string:}{}
    \SetKwFunction{Range}{range}
    \SetKw{KwTo}{in}\SetKwFor{For}{for}{\string:}{}%
    \SetKwIF{If}{ElseIf}{Else}{if}{:}{elif}{else:}{}%
    \SetKwFor{While}{while}{:}{fintq}%
    \AlgoDontDisplayBlockMarkers\SetAlgoNoEnd\SetAlgoNoLine%
    $\Vx^{k+\frac{1}{2}} = \Vx^{k+1} - \alpha_{f} \Vp$\;
    $\Vg = \Vg - \alpha_{f} \Vp$\;
    $\Vx^{k+1} = P_{\Omega} (\Vx^{k+\frac{1}{2}} -\overline{\alpha} \tilde{\bm{d}})$\;
    $\Vg = \MA \Vx^{k+1} -\Vb$\;
    $\Vp = \Vphi$
    \caption{Expansion}
    \label{alg:exp}
\end{algorithm}

\begin{algorithm}[htb]
    \SetStartEndCondition{ }{}{}%
    \DontPrintSemicolon
    \SetKwProg{Fn}{def}{\string:}{}
    \SetKwFunction{Range}{range}
    \SetKw{KwTo}{in}\SetKwFor{For}{for}{\string:}{}%
    \SetKwIF{If}{ElseIf}{Else}{if}{:}{elif}{else:}{}%
    \SetKwFor{While}{while}{:}{fintq}%
    \AlgoDontDisplayBlockMarkers\SetAlgoNoEnd\SetAlgoNoLine%
    $\alpha_{cg} = \Vg^T \Vbeta / (\Vbeta)^T \MA \Vbeta$\;
    $\Vx^{k+1} = \Vx^{k} -\alpha_{cg}\Vbeta$\;
    $\Vg = \Vg - \alpha_{cg} \MA \Vbeta$\;
    $\Vp = \Vphi$\;
    \caption{Proportioning}
    \label{alg:prop}
\end{algorithm}

Let us briefly explain the algorithm. In each iteration, the algorithm checks that the current approximation of solution $\Vx^k$ is strictly proportional
\begin{equation}
\label{eq:prop}
    ||\Vbeta(\Vx^k)||^2 \le \Gamma^2 ||\Vg^f(\Vx^k)||^2, \qquad \Gamma > 0.
\end{equation}
If this inequality does not hold, the chopped gradient $\Vbeta$ dominates (depending on the value of $\Gamma$, typically $\Gamma=1$) the norm of the projected gradient $\Vg^P$ and therefore we need to release some components from the active set by a proportioning step.
The proportioning step consists of a single steepest descent step in the direction opposite $\Vbeta$.

On the other hand, if the current solution is proportional, i.e., \eqref{eq:prop} holds, then the free gradient $\Vphi$ dominates the norm of $\Vg^P$, and we focus on minimization of $\Vphi$.
First, we compute $\alpha_{cg}$ as the optimal step length for minimization in direction $-\Vp$ and $\alpha_f$ as the maximal step length in this direction that does not leave the feasible set.
If $\alpha_{cg} \le \alpha_f$ we can do an unconstrained minimization using a standard CG step; otherwise, we do the expansion step.
Note that initially and after both expansion and proportioning steps $\Vp = \Vphi$, while the CG steps set the next minimization direction A-orthogonal to the previous one.

The expansion consists of a so-called half-step which is the step with maximal step length $\alpha_f$.
The half-step expands the active set, typically, by one component.
Then a step in the direction opposite $\tilde{\mb{d}}$, where $\tilde{\mb{d}} = \Vg^f$ or $\tilde{\mb{d}} = \Vg^r$, with a fixed step length $\overline{\alpha} \in (0,2||\MA||^{-1}]$ is performed.
Notice that due to the construction of $\Vg^r$, we have
\begin{displaymath}
  \Vx^{k+1} = P_{\Omega} (\Vx^{k+\frac{1}{2}} -\overline{\alpha} \Vphi) = \Vx^{k+\frac{1}{2}} -\overline{\alpha} \bm{g}^r
\end{displaymath}
The active set is expanded in a component $j$ if this component is in the free set and
\begin{displaymath}
  g_j > 0 \quad \text{and} \quad \overline{\alpha} g_j \ge x_j - l_j
\end{displaymath}
or
\begin{displaymath}
  g_j \le 0 \quad \text{and} \quad \overline{\alpha} g_j \le x_j - u_j.
\end{displaymath}
Therefore, $\overline{\alpha}$ controls how large a component of gradient (in the correct direction) has to be to expand the active set in the given component.
Larger values of $\overline{\alpha}$ can potentially expand the active set in a greater number of components.
However, even with the largest possible value of $\overline{\alpha}$, the active set may not be expanded at all.
As will be demonstrated in the next section, the expansion step also decreases the cost function.

The operation count for each of the three steps is summarised in Table \ref{tab:opcnt}. We would argue that in most cases, the cost of a step primarily depends on the number of Hessian multiplication it performs.

Note that, either bound can be omitted in the formulation of the algorithm. If both bounds are omitted, the algorithm is equivalent to a standard CG method.

\begin{table}
\centering
\begin{tabular}{|l|c|c|c|c|c|}
\hline
Step & Hess. mult. & Dot prod. &  Vec. update & Grad. split.\\\hline
CG & 1 & 2 & 3 & 1 \\\hline
Expansion & 2 & 1 & 5 & 2\\\hline
Proportioning & 1 & 1 & 3 & 1 \\\hline\hline
Expansion-optapprox & 2 & 3 & 5 & 2\\\hline
Expansion-opt & 3 & 3 & 5 & 2\\\hline
Expansion-projCG & 2 & 1 & 3 & 1 \\ \hline
\end{tabular}
\caption{Number of operations per MPRGP step. The bottom half of the table contains the newly proposed variants of the expansions step.}
\label{tab:opcnt}
\end{table}

\section{Expansion Modifications}
\label{sec:expansion}
As we discussed in the previous section, the time to solution is primarily determined by the number of Hessian multiplications. Therefore, we need to minimize the overall number of Hessian multiplications to speed-up the MPRGP algorithm. 

Numerical experiments (e.g., in Section \ref{sec:num}) show that MPRGP may need many expansion steps to identify the active set because standard expansion steps often enlarge the active set by only one or a few components.
Moreover, expansion is about twice as expensive as the other steps.
We can try modifying the expansion step to, potentially, enlarge the active set faster.
Such modifications should lead to a decrease in the number of expansion steps as well as the overall number of Hessian multiplications.

\subsection{Expansion Step with Adaptive Step Length}


As our experiments demonstrate, it is better to perform longer steps which are often unfeasible and to make a subsequent projection onto the feasible set.
One of the ideas how to choose these longer steps comes from the steepest descent method.
We analyze the step lengths in a given direction to find a step length for which the cost function decrease is maximal.
As our cost function is parabolic, we can guarantee the cost function decrease with up to twice the optimal step length.

Let us reiterate that our cost function is 
$$
  f(\mb{x})=\frac{1}{2}\mb{x}^{T}\mb{Ax}-\mb{x}^{T}\mb{b},
$$
and that the expansion does a step in the $\mb{g}^r$ direction with a fixed step length $\overline{\alpha} \in (0,2||\MA||^{-1}]$.
Let us assume that we do a step in a direction $\mb{d}$ instead, where $\mb{d}$ is either $\mb{g}^r$ or $\mb{g}^f$ so that no active component is freed.
We want to choose a step length such that the cost function decreases, i.e.,
\begin{eqnarray*}
\begin{array}{lcl}
  f(\mb{x}) - f(\mb{x}-\bar{\alpha}\mb{d}) & = & f(\mb{x}) - \frac{1}{2}(\mb{x}-\bar{\alpha}\mb{d})^{T}\mb{A}(\mb{x}-\bar{\alpha}\mb{d})+(\mb{x}-\bar{\alpha}\mb{d})^{T}\mb{b}= \\[2mm]
  & = & f(\mb{x}) - \frac{1}{2}\mb{x}^{T}\mb{Ax}+\mb{x}^{T}\mb{b}-\frac{1}{2}\bar{\alpha}^{2}\mb{d}^{T}\mb{A}\mb{d}+\bar{\alpha}\mb{d}^{T}\mb{Ax}-\bar{\alpha}\mb{d}^{T}\mb{b}=\\[2mm]
  & = & -\frac{1}{2}\bar{\alpha}^{2}\mb{d}^{T}\mb{A}\mb{d}+\bar{\alpha}\mb{d}^{T}\mb{g}\ge 0.
\end{array}
\end{eqnarray*}
and after division by $\bar{\alpha} > 0$
$$
  \frac{1}{2}\bar{\alpha}\mb{d}^{T}\mb{A}\mb{d} \leq \mb{d}^{T}\mb{g}.
$$

Assuming $\mb{d}$ is not in the null space of $\MA$, we have  $\mb{d}^{T}\mb{Ad} > 0$ and so we can divide the inequality by $\mb{d}^{T}\mb{A}\mb{d}$
$$
  \bar{\alpha} \leq \frac{2\mb{d}^{T}\mb{g}}{\mb{d}^{T}\mb{A}\mb{d}}.
$$

Because $\mb{d}$ is either free gradient or reduced free gradient, we have 
$$
  \mb{d}^{T}\mb{d} \leq \mb{d}^{T}\mb{g} \leq \mb{g}^{T}\mb{g}
$$
so that 
$$
\bar{\alpha} \leq \frac{2\mb{d}^{T}\mb{d}}{\mb{d}^{T}\mb{A}\mb{d}} \leq \frac{2\mb{d}^{T}\mb{g}}{\mb{d}^{T}\mb{A}\mb{d}}.
$$

We know that 
$$
  \mb{d}^{T}\mb{d} = \left\| \mb{d}\right\|^{2}\text{{\,\,\,\,\, and\,\,\,\,\,}}\mb{d}^{T}\mb{A}\mb{d} = \left|\mb{d}^{T}\mb{A}\mb{d}\right| \leq \left\| \mb{A}\right\| \left\| \mb{d}\right\|^{2}
$$
that proves 
$$
  \frac{2\mb{d}^{T}\mb{d}}{\mb{d}^{T}\mb{A}\mb{d}} \geq \frac{2\left\| \mb{d}\right\|^{2}}{\left\| \mb{A}\right\| \left\| \mb{d}\right\|^{2}}=2\left\| \mb{A}\right\|^{-1}.
$$

Furthermore 
$$
  \bar{\alpha} \leq \frac{2\mb{d}^{T}\mb{g}}{\mb{d}^{T}\mb{A}\mb{d}} = \frac{2\mb{d}^{T}\mb{g}}{\mb{d}^{T}\mb{A}\mb{d}} \cdot 1 = \frac{2\mb{d}^{T}\mb{g}}{\mb{d}^{T}\mb{A}\mb{d}} \cdot \frac{\mb{d}^{T}\mb{d}}{\mb{d}^{T}\mb{d}} = \frac{2\mb{d}^{T}\mb{d}}{\mb{d}^{T}\mb{A}\mb{d}} \cdot \frac{\mb{d}^{T}\mb{g}}{\mb{d}^{T}\mb{d}},
$$
$$
  \frac{\mb{d}^{T}\mb{g}}{\mb{d}^{T}\mb{d}} \geq 1
$$
which gives us a larger upper bound for $\bar{\alpha}$
$$
  0 \leq \bar{\alpha} \leq 2\left\| \mb{A}\right\|^{-1} \leq \frac{2\mb{d}^{T}\mb{d}}{\mb{d}^{T}\mb{A}\mb{d}} \leq 2\left\| \mb{A}\right\|^{-1}\frac{\mb{d}^{T}\mb{g}}{\mb{d}^{T}\mb{d}} \leq \frac{2\mb{d}^{T}\mb{g}}{\mb{d}^{T}\mb{A}\mb{d}}.
$$

We can consider this bound 
$$
  0 \leq \bar{\alpha} \leq 2\left\| \mb{A}\right\|^{-1}\frac{\mb{d}^{T}\mb{g}}{\mb{d}^{T}\mb{d}} \leq \frac{2\mb{d}^{T}\mb{g}}{\mb{d}^{T}\mb{A}\mb{d}}
$$
as an adaptive step length for an expansion step taking into account the actual situation.
For testing, let us consider the following notation:
\begin{itemize}
  \item \emph{fixed} $\bar{\alpha} = \alpha_u\left\| \mb{A}\right\|^{-1},$
  \item \emph{optapprox} $\bar{\alpha} = \alpha_u\left\| \mb{A}\right\|^{-1}\frac{\mb{d}^{T}\mb{g}}{\mb{d}^{T}\mb{d}},$
  \item \emph{opt} $\bar{\alpha} = \alpha_u \frac{\mb{d}^{T}\mb{g}}{\mb{d}^{T}\mb{A}\mb{d}}$,
\end{itemize}
where $\alpha_u \in (0,2]$.
Any choice of the $\overline{\alpha}$ from the above options guarantees the reduction of the cost function.
However, only the fixed choice of $\overline{\alpha}$ and descent direction $\mb{d} = \mb{g}^r$ ensures that $\mb{x}$ is kept in the feasible set.
Therefore, we need to project the new approximation to the feasible set, i.e,
\begin{displaymath}
  \Vx^{k+1} = P(\Vx^{k+\frac{1}{2}} -\overline{\alpha}\tilde{\mb{d}}).
\end{displaymath}
The above expression is exactly the expansion step of the original algorithm with $\tilde{\mb{d}} = \Vg^f$. 
Since we derived the step lengths without taking into considerations the projection, it makes sense to decouple the $\mb{d}$ used for computation of $\overline{\alpha}$ and the descent direction $\tilde{\mb{d}}$.

Note that some combinations of vectors used for step length computation as well as descent directions are equivalent.
Namely, \emph{fixed} step length is equivalent to \emph{optapprox} with $\mb{d} = \Vg^f$ with any of the two descent directions.

\subsection{Expansion Using Projected CG Step}
Recall, that the expansion consists of the half step followed by the expansion step line search.
The half step is a CG step with step length reduced such that the computed approximation is in the feasible set.

Since our goal is to expand the active set faster, it seems reasonable to replace the half step by the full CG step with a subsequent projection onto the feasible set.
To be more specific, our expansion step becomes
\begin{displaymath}
  \Vx^{k+1} = P_{\Omega}(\Vx^k -\alpha_{cg} \Vp),
\end{displaymath}
followed by reseting $\Vp = \Vg^f$.
Note, that realising the expansion in this way simplifies the implementation as we can always compute the CG step and then compute the gradient using CG recurrence when the step was feasible; otherwise, we project the solution onto the feasible set and recompute the gradient explicitly. Algorithm \ref{alg:projcg} illustrates the implementation. It replaces \emph{if\dots else} block on lines 6--9 in Algorithm \ref{alg:mprgp}.

\begin{algorithm}[htb]
    \SetStartEndCondition{ }{}{}%
    \DontPrintSemicolon
    \SetKwProg{Fn}{def}{\string:}{}
    \SetKwFunction{Range}{range}
    \SetKw{KwTo}{in}\SetKwFor{For}{for}{\string:}{}%
    \SetKwIF{If}{ElseIf}{Else}{if}{:}{elif}{else:}{}%
    \SetKwFor{While}{while}{:}{fintq}%
    \AlgoDontDisplayBlockMarkers\SetAlgoNoEnd\SetAlgoNoLine%
    $\Vx^{k+1} = \Vx^{k} -\alpha_{cg}\Vp$\;
    \If{$\alpha_{cg} \le \alpha_{f}$}{
      $\Vg = \Vg - \alpha_{cg} \MA \Vp$\;
      $\beta = \Vp^T \MA \Vphi / \Vp^T \MA \Vp$\;
      $\Vp = \Vphi -\beta \Vp$\;
    }\Else{
      $\Vx^{k+1} = P_{\Omega}(\Vx^{k+1})$\;
      $\Vg = \MA \Vx^{k+1} -\Vb$\;
      $\Vp = \Vphi$\;
    }
    \caption{Projected CG}
    \label{alg:projcg}
\end{algorithm}

\begin{figure}
  \centering
  \includegraphics[width=.5\textwidth]{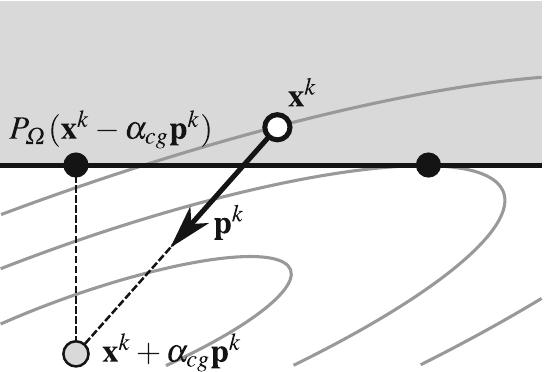}
  \caption{Projected CG step illustration after the first unconstrained CG step. The grey part is the feasible set bounded by the horizontal line. The ellipses are the cost function contour lines. The right black point on the horizontal line is the solution. \cite{Dos-book-17}}
    \label{fig:projcg}
\end{figure}

However, in some cases, this step can lead to an increase in the cost function.
In \cite{Dos-book-17}, section 6.2, the author illustrates by Figure \ref{fig:projcg} the situation when increase in the cost function happens.
Clearly, the CG step finds the uncostrained minimizer of the cost function, but the subsequent projection onto the feasible set puts as on a higher contour line, i.e., increases the value of cost function.
In this case the standard expansion would put as closer to solution.
Yet, unless $\overline{\alpha}$ happens to be such that the line search finds the exact solution, both approaches converge in the next iteration.
Note that in this example, \emph{opt} step length with $\alpha_u = 1$ with any combination of the allowed vector for both the compuation of the step length and line search direction would converge to the exact solution in a single iteration.

\section{Short Introduction of Software and Methods}
\label{sec:soft}
We used PERMON \cite{permon-www,HapHorPos-HPCSE-16} for the numerical experiments. PERMON is a collection of open-source software libraries used for quadratic programming (QP) and its applications. It is based on PETSc \cite{petsc-web-page,petsc-user-ref} and follows the same design and coding style, making it easy to use for anyone familiar with PETSc. 

The main module is PermonQP. It provides data structures, transformations, solvers, and supporting functions for QP. The transformations, e.g. dualization, can be used to simplify a QP problem and provides functions to reconstruct the solution of the original problem. Among the solvers available are augmented Lagrangian-based algorithms (e.g., SMALE \cite{Dos-book-09}), gradient projection-type methods (e.g., MPRGP), solvers available in TAO \cite{tao-user-ref}, and others.

PermonFLLOP (FETI Light Layer on Top of PETSc) implements domain decomposition methods of the FETI type. It applies PermonQP transformations on a primal problem to derive either unconstrained or box and equality constrained FETI-type QP formulation for unconstrained and contact problems, respectively. PermonQP then solves the final formulation. Moreover, it includes functions for generating subdomain glueing matrix, identifying subdomain kernel, and efficient coarse problem solution. 

The PermonSVM \cite{permonsvm-www,Kruzik-SVM-LNCS-2018} package provides an implementation of binary classification via soft-margin
Support Vector Machines (SVMs). It implements a scalable training procedure based on a linear kernel, taking advantage of an implicit representation of the Gramm matrix. It utilises PermonQP to solve the dual SVM formulation.

\subsection{TFETI for contact problems}

Let us consider the a spatial domain $\Omega$ which is decomposed into non-overlapping subdomains.
Then virtually arbitrary Finite Element Method (FEM) implementation can be used to generate the subdomain stiffness matrices $\mb K_s$ and the subdomain load vectors $\mb f_s$
as sequential data for each subdomain $\Omega_s$, $s=1,\,\ldots,\,N_S$ independently.

The original primal problem 
\begin{equation}
  \argmin_{\mb{u}}\frac{1}{2}\mb{u}^{T}\mb{Ku}-\mb{f}^{T}\mb{u}\;\;\text{s.t.}\;\;\mb{B}_{I}\mb{u}\leq\mb{o}\;\;\text{and}\;\;\mb{B}_{E}\mb{u}=\mb{o},\label{ekv_16-1}
\end{equation}
where $\mb K = diag(\mb{K}_1,\dots,\mb{K}_{N_S})$ is global stiffness matrix, $\mb{f} =\left[\mb{f}_1^T, \dots ,\mb{f}_{N_S}^T\right]^T$ is global right hand side, $\mb u$ is unknown displacement, $\mb B_I$ represents non-penetration condition, and $\mb B_E$ glues the subdomains together. The primal problem is transformed into dual one
\begin{equation}
  \argmin_{\mb{u}}\frac{1}{2}\mb{\lambda}^{T}\mb{F\lambda}-\mb{\lambda}^{T}\mb{d}\;\;\text{s.t.}\;\;\mb{\lambda}_{I}\geq\mb o\:\;\text{and}\;\:\mb{G\lambda=e},\label{ekv_24}
\end{equation}
where
$$
\mb{G}=\mb{R}^{T}\mb{B}^{T},\,\mb{d}=\mb{B}\mb{K}^{\dagger}\mb{f},\,\mb{e}=\mb{R}^{T}\mb{f},\,\mb{F}=\mb{B}\mb{K}^{\dagger}\mb{B}^{T},
$$
$\mb{K}^{\dagger}$ denotes a left generalized
inverse of $\mb{K}$, i.e. a matrix satisfying $\mb{K}\mb{K}^{\dagger}\mb{K}=\mb{K}$ and columns of $\mb R$ span the null space.
The constraint matrix 
$\mb{B}=\left[\begin{smallmatrix}\mb{B}_{I}^T & \mb{B}_{E}^T\end{smallmatrix}\right]^T$
can be constructed so that it has a full rank, and then the Hessian $\mb F$ is positive definite with a relatively
favourably distributed spectrum for application of the CG method. For more details see ...

\subsection{SVM and no-bias data classifications}

Support Vector Machines (SVMs) belong to the conventional machine learning (ML) techniques, and they can solve classification as well as regression problems.
Despite the fact that deep learning (DL) \cite{Lecun-Nature-2015} is getting popular in the recent years, SVMs are successfully applied for specific tasks in various scientific areas including genetics \cite{Brown-NASUSA-2000}, geosciences \cite{Shi-IEEE-2012}, and image analysis \cite{Foody-RSE-2006}.
Unlike the DL underlying architecture, SVMs could be considered as the single perceptron problems that find the learning functions that maximize the geometric margins.
Therefore, we can simply explain the qualities of a learning model and the underlying solver behaviour. 
In this paper, we will focus on the linear SVMs for classifications.

SVM was originally designed as a supervised binary classifier \cite{CorVap-ML-1995},  i.e., a classifier that decides whether a sample falls into either Class A or Class B employing a model determined from already categorised samples in the training phase of the classifier.
Let us denote the training data as an ordered sample-label pairs such that 
\[ T := \{ \left(\bm{x}_1, y_1\right), \ \left(\bm{x}_2, y_2\right), \ \dots, \ \left(\bm{x}_m,  y_m\right) \},\]
where $m$ is the number of samples, $\bm{x}_i \in \mathbb{R}^n, \ n \in \mathbb{N},$ is the $i$-th sample and $y_i \in \{-1, 1\}$ denotes the label of the $i$-th sample, $i \in \{1, \ 2, \ \dots, m\}$.
The label determines the sample's class.

Standard SVM solves problem of finding a classification model in a form of maximal-margin hyperplane  $H = \left\langle\bm{w}, \bm{x}\right\rangle + b$, where $\bm{w}$ is the normal vector of the hyperplane $H$ and $b$ is its bias from origin.
In the no-bias classification, we do not consider  bias $b$ in a classification model, but we include it into the problem by means of augmenting the vector $\bm{w}$ and each sample $\bm{x}_i$ with an additional dimension so that $\bm{\widehat{w}} \leftarrow \begin{bmatrix}\bm{w} \\ b \end{bmatrix}$, $\bm{\widehat{x}_i} \leftarrow \begin{bmatrix}\bm{x}_i \\ \beta \end{bmatrix}$, where $\beta \in \mathbb{R}^+$ is a user defined variable (typically set to $1$).

The problem of finding hyperplane $\widehat{H} = \left\langle \bm{\widehat{w}}, \bm{\widehat{x}} \right\rangle$ can be formulated as a constrained optimization problem in the following primal formulation
\begin{equation}
\label{eq:svmPrimal}
	 \ \argmin_{\bm{\widehat{w}}, \ \xi_i} \ \frac{1}{2} \left\langle \bm{\widehat{w}}, \bm{\widehat{w}} \right\rangle + \frac{C}{p} \sum_{i = 1}^n \xi_i^p \ \text{s.t.} \
		  \begin{cases}
			 \ y_i\left\langle\bm{\widehat{w}}, \bm{\widehat{x}}_i\right\rangle \geq 1 - \xi_i, \ i \in \{1,2, \dots, n\}, \\
			 \ \xi_i \geq 0, \  i \in \{1,2, \dots, n\},
		  \end{cases}
\end{equation}
where  $p \in \{1, 2\}$, $\xi_i =\max\left(0, 1 - y_i\left\langle\bm{\widehat{w}}, \bm{\widehat{x}}_i\right\rangle\right)$ is the hinge loss function and $C \in \mathbb{R}^+$ is a user defined penalty. Using the Lagrange duality and denoting $\bm{H} = \bm{Y}^T \bm{G} \bm{Y}$, $\mb{Y} = diag(\mb{y})$, $\mb{y} = \left[y_1, \ y_2, \  \dots, y_m\right]^T$, $\bm{G} = \bm{X}^T\bm{X}$, $\mb{X}= \begin{bmatrix} \mb{\widehat{x}}_1 & \dots & \mb{\widehat{x}}_m \end{bmatrix}$, $\bm{e} = \left[1, 1, \dots, 1\right]^T \in \mathbb{R}^m$, $\bm{o} = \left[0, \ 0,  \ \dots, \ 0\right]^T \in \mathbb{R}^m$ 
, we transform \eqref{eq:svmPrimal} into the dual formulations  
\begin{equation}
\label{eq:l1svm}
	\argmin_{\boldsymbol{\lambda}} \ \frac{1}{2} \boldsymbol{\lambda}^T \bm{H}  \boldsymbol{\lambda} -  \boldsymbol{\lambda}^T \bm{e}  \quad \text{s.t.} 
	\quad	\boldsymbol{o} \leq  \boldsymbol{\lambda} \leq C\bm{e}
\end{equation}
for $p = 1$ and 
\begin{equation}
\label{eq:l2svm}
     \argmin_{\boldsymbol{\lambda}} \ \frac{1}{2} \boldsymbol{\lambda}^T \left(\bm{H} + C^{-1} \bm{I} \right) \boldsymbol{\lambda} -  \boldsymbol{\lambda}^T \bm{e}  \quad \text{s.t.} 
	\quad	\boldsymbol{o} \leq  \boldsymbol{\lambda}
\end{equation}	
for $p = 2$.
The first formulation is commonly called dual $l1$-loss SVM, and the second one is known as dual $l2$-loss SVM.
The Hessian associated with QP  problem \eqref{eq:l1svm} is symmetric positive-semi definite. It becomes positive definite using regularization by matrix $C^{-1}\bm{I}$ in formulation \eqref{eq:l2svm}.

\section{Numerical experiments}
\label{sec:num}
This section compares the presented expansion variants on two benchmarks.
The first one is a 3D linear elasticity contact problem using TFETI.
The second benchmark consists of several classification problems solved by SVMs.

In each benchmark, MPRGP parameter $\Gamma = 1$.

\subsection{3D Linear Elasticity Contact Problem}
\begin{figure}
    \includegraphics[width=\textwidth]{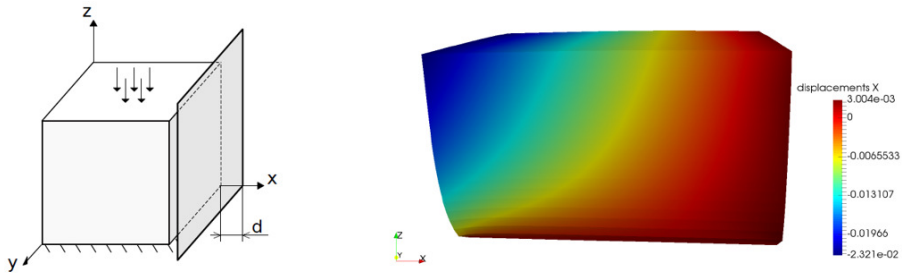}
    \caption{3D elastic cube contact problem.}
    \label{fig:cube}
\end{figure}

The first benchmark is a 3D linear elasticity contact problem.
We considered an elastic cube with the dimensions $1\times 1\times 1$ [mm] with the bottom face fixed, the top one loaded with a vertical surface force $f_{z}=-465$ [N/mm$^{2}$] directed downwards, and the right one in contact with a rigid obstacle. Young modulus is $E=2\cdot10^{5}$ [MPa] and Poisson ratio $\mu=0.33$. See Figure \ref{fig:cube}.

We use TFETI domain decomposition to solve this problem.
Let us consider only one regular decomposition into 1000 subdomains (10 in each direction) with 27,000 elements per subdomains (30 in each direction) making a total of 81,812,703  (undecomposed) degrees of freedom.
We use M-variant of the SMALBE \cite{Dos-book-09} algorithm to take care of the equality constraint with MPRGP used as the inner solver.
The stopping tolerance of the outer solver (SMALBE) is set to $10^{-6}$ relative to the right-hand side.
Other parameters for SMALBE are $M = 100||\MA||$, $\eta = 1.1||\MA||$ and, when the decrease of the Lagrangian is sufficient, $M$ is reduced by $10$.

\begin{figure}
    \includegraphics[width=\textwidth]{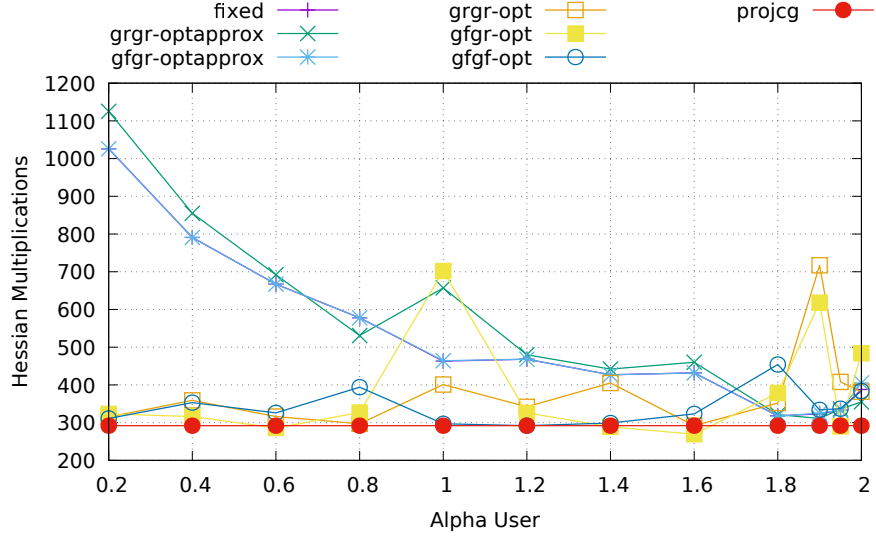}
    \caption{Linear elasticity contact problem: Comparison of expansion strategies in the term of the number of the Hessian multiplications depending on $\alpha_u$.}
    \label{fig:gcube}
\end{figure}

Results, in term of the overall number of Hessian multiplications are reported in Figure \ref{fig:gcube}. Note, that the projected CG expansion variant does not contain parametr $\alpha_u$. 
Table \ref{tab:cube} compares the numbers of outer and inner iterations, i.e. numbers of CG, expansion, proportioning steps and Hessian multiplications for the three best values of $\alpha_{u}$. The full results are reported in \ref{sec:appendix}.

The standard expansion with the fixed-length expansion step needs $318$ Hessian multiplication for the best $\alpha_u = 1.8$.
The number of Hessian multiplications mostly decreases until the best $\alpha_u$ is reached and then starts to increase again.

The \emph{optapprox} variants do not bring much as they achieve the convergence in a very similar number of iterations to the \emph{fixed} variant.
This behaviour suggests that the free gradient is not sufficiently different from the reduced free gradient, i.e., $\Vg^f \approx \Vg^r$ and we have
\begin{displaymath}
  \overline{\alpha} = \alpha_u ||\MA||^{-1} \frac{(\Vg^r)^T\Vg}{(\Vg^r)^T\Vg^f} \approx \alpha_u ||\MA||^{-1} \frac{(\Vg^f)^T\Vg}{(\Vg^f)^T\Vg^f} = \alpha_u ||\MA||^{-1}.
\end{displaymath}

Much better are the \emph{opt} step lengths.
They can significantly decrease the number of expansion steps needed (up to half in one case).
The number of Hessian multiplications is decreased as well, but due to the additional Hessian multiplication, the effect is less pronounced.
The best result of $269$ Hessian multiplication achieved step length computed with $\Vg^r$ in the direction of $\Vg^f$ for $\alpha_u=1.6$, a decrease of $15 \%$ compared to the best \emph{fixed} variant).
The other variants of \emph{opt} outperformed \emph{fixed} as well, requiring $292$ Hessian multiplication (reduction of $8 \%)$.
The problem with any of the \emph{opt} strategies is that there are quite large jumps in the number of Hessian multiplications depending on the value of $\alpha_u$. 

Lastly, the project CG step variant of expansion performs well. It achieves $292$ Hessian multiplication, a reduction of $8 \%$ compared to \emph{fixed}. Moreover, it is outperformed only by the best \emph{opt} strategy, but there is no dependence on the $\alpha_u$. 


\begin{table}
\centering
\begin{tabular}{| *{7}{c|}}
\hline
exp. type & $\alpha_u$ & outer it. & \#Hess. mult. & \#CG  & \#Exp.  & \#Prop.  \\ 
\hline
\multirow{3}{*}{fixed}
 & 1.8 & 10 & 318 & 158 & 74 & 2 \\
 & 1.9 & 10 & 323 & 144 & 83 & 3 \\
 & 1.95 & 7 & 333 & 70  & 127 & 2 \\
\hline
 \multirow{3}{*}{grgr-optapprox}
 & 1.9 & 10 & 311 & 162 & 68 & 3 \\
 & 1.8 & 10 & 321 & 145 & 82 & 2 \\
 & 1.95 & 10 & 336 & 125 & 99 & 3 \\
\hline
\multirow{3}{*}{gfgr-optapprox}
 & 1.8 & 10 & 317 & 159 & 73 & 2 \\
 & 1.9 & 10 & 325 & 140 & 86 & 3 \\
 & 1.95 & 8 & 333 & 95 & 114 & 3 \\
\hline
 \multirow{3}{*}{grgr-opt}
 & 1.6 & 8 & 292 & 120 & 54 & 2 \\
 & 0.8 & 9 & 297 & 132 & 51 & 3 \\
 & 0.2 & 8 & 315 & 113 & 64 & 2 \\
\hline
 \multirow{3}{*}{gfgr-opt}
 & 1.6 & 9 & 269 & 129 & 43 & 2 \\
 & 0.6 & 12 & 286 & 139 & 44 & 3 \\
 & 1.4 & 9 & 289 & 154 & 41 & 3 \\
\hline
 \multirow{3}{*}{gfgf-opt}
 & 1.2 & 10 & 292 & 151 & 43 & 2 \\
 & 1.0 & 11 & 297 & 171 & 37 & 4 \\
 & 1.4 & 9 & 299 & 152 & 45 & 3 \\
\hline
 projcg & - & 10 & 292 & 171 & 53 & 5 \\
\hline
\end{tabular}
\caption{Linear elasticity contact problem: Comparison of expansion strategies by the number of SMALBE outer iterations, overall Hessian multiplications, CG, Expansion, and proportioning steps.}
\label{tab:cube}
\end{table}

\subsection{Classification Problems}

The second benchmark uses SVMs for classifications on three publicly available datasets, namely Australian, Diabetes, and Ionosphere downloaded from LIBSVM dataset webpage \cite{www-libsvm-datasets}.
The Australian dataset (Australian Credit Approval) concerns credit card applications. The objective of the Diabetes dataset is to predict whether a patient has diabetes.
Finally, the Ionosphere dataset is about the classification of radar returns from the ionosphere as either suitable for further analysis or not. 
The information about the number of samples and features for each dataset are in Table \ref{tab:svm.datasets}. 

\begin{table}
\centering
\begin{tabular}{|l|r|r|}
\hline
Dataset & \# samples & \# features \\ 
\hline
Australian &  690 &  14 \\ \hline
Diabetes   &  678 &   8 \\ \hline
Ionosphere &  351 &  34 \\ \hline
\end{tabular}
\caption{The number of samples and features for our classification problem datasets.}
\label{tab:svm.datasets}
\end{table}

We use the no-bias SVM formulation with $l1$ hinge-loss fuction and with $C=1$.
Since this is a classification application, there are very low requirements on the accuracy of the underlying solver \cite{Kruzik-SVM-LNCS-2018}.
Therefore, the relative tolerance of MPRGP is set to $10^{-1}$.
The initial guess is set just under the upper bound (each component is set to $1-100\epsilon_m$, where $\epsilon_m \approx 2.2\mathrm{e}{-16}$ is the machine epsilon).

\begin{figure}
    \includegraphics[width=\textwidth]{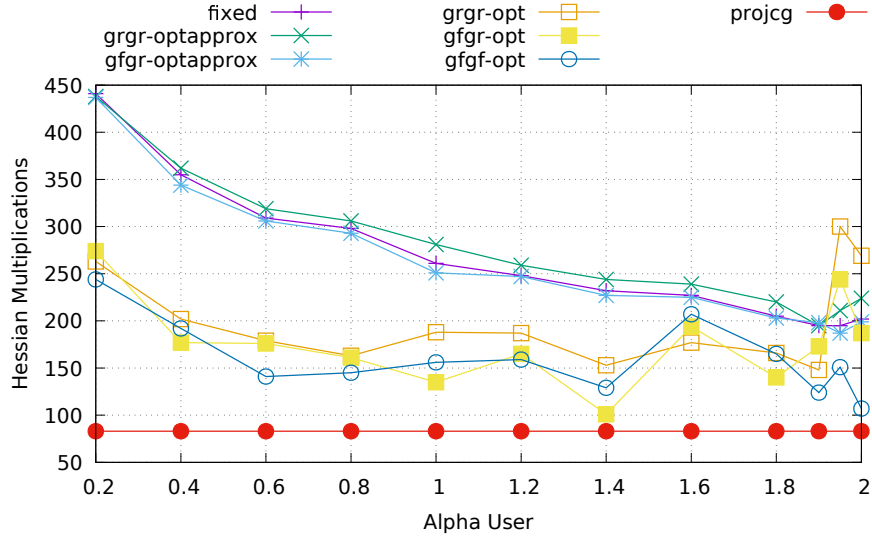}
    \caption{Classification problem, Australian dataset: Comparison of expansion strategies in the term of the number of the Hessian multiplications depending on $\alpha_u$.}
    \label{fig:gaustralian}
\end{figure}

\begin{figure}
    \includegraphics[width=\textwidth]{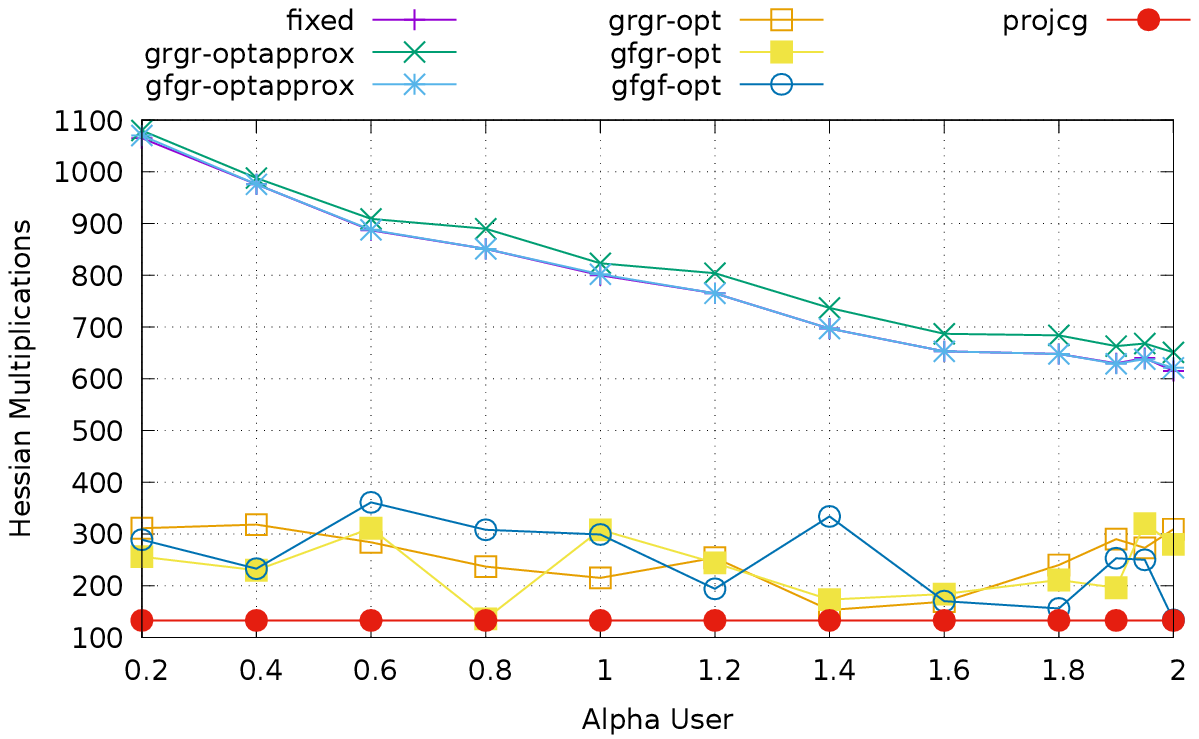}
    \caption{Classification problem, Diabetes dataset: Comparison of expansion strategies in the term of the number of the Hessian multiplications depending on $\alpha_u$.}
    \label{fig:gdiabetes}
\end{figure}

\begin{figure}
    \includegraphics[width=\textwidth]{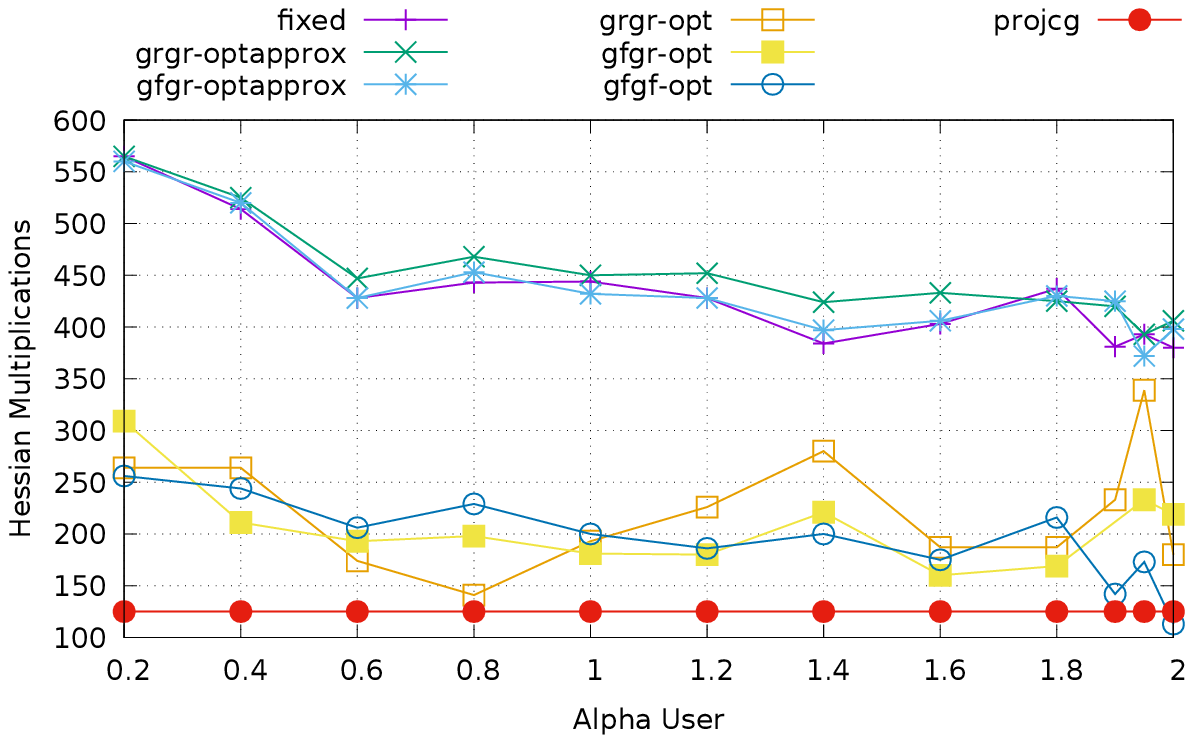}
    \caption{Classification problem, Ionosphere dataset: Comparison of expansion strategies in the term of the number of the Hessian multiplications depending on $\alpha_u$.}
    \label{fig:gionosphere}
\end{figure}

As in the previous bechmark, we report the results in term of the overall number of Hessian multiplications in Figures \ref{fig:gaustralian}, \ref{fig:gdiabetes}, and \ref{fig:gionosphere}.
Note that we removed an outlier from the Ionosphere dataset graph; for \emph{gfgr-opt} with $\alpha_u = 1.9$ there was a total of $2,322$ Hessian multiplications required for convergence.
In Tables \ref{tab:australian}, \ref{tab:diabetes}, and \ref{tab:ionosphere} the number of CG, expansion, proportioning steps and Hessian multiplications are compared for the three best values of $\alpha_{u} = 1.9$.  The full results are reported in \ref{sec:appendix}.

\begin{table}
\centering
\begin{tabular}{| *{6}{c|}}
\hline
exp. type & $\alpha_u$ & \#Hess. mult. & \#CG  & \#Exp.  & \#Prop.  \\ 
\hline
\multirow{3}{*}{fixed}
 & 1.9  & 195 & 8 & 92 & 2 \\
 & 1.95 & 195 & 8 & 92 & 2 \\
 & 2.0  & 202 & 10 & 95 & 1 \\
\hline
 \multirow{3}{*}{grgr-optapprox}
 & 1.9  & 195 & 9 & 92 & 1 \\
 & 1.95 & 211 & 7 & 101 & 2 \\
 & 1.8  & 220 & 14 & 102 & 1 \\
\hline
\multirow{3}{*}{gfgr-optapprox}
 & 1.95 & 187 & 5  & 90 & 1 \\
 & 1.9  & 198 & 16 & 90 & 1 \\
 & 2.0  & 199 & 5  & 96 & 1 \\
\hline
 \multirow{3}{*}{grgr-opt}
 & 1.9 & 148 & 28 & 39 & 2 \\
 & 1.4 & 153 & 27 & 41 & 2 \\
 & 0.8 & 163 & 16 & 48 & 2 \\
\hline
 \multirow{3}{*}{gfgr-opt}
 & 1.4 & 101 & 20 & 26 & 2 \\
 & 1.0 & 135 & 15 & 39 & 2 \\
 & 1.8 & 140 & 25 & 37 & 3 \\
\hline
 \multirow{3}{*}{gfgf-opt}
 & 2.0 & 107 & 23 & 27 & 2 \\
 & 1.9 & 124 & 28 & 31 & 2 \\
 & 1.4 & 129 & 28 & 33 & 1 \\
\hline
 projcg & - & 83 & 16 & 32 & 2 \\
\hline
\end{tabular}
\caption{Classification problem, Australian dataset: Comparison of expansion strategies by the number of SMALBE outer iterations, overall Hessian multiplications, CG, Expansion, and proportioning steps.}
\label{tab:australian}
\end{table}

\begin{table}
\centering
\begin{tabular}{| *{6}{c|}}
\hline
exp. type & $\alpha_u$ & \#Hess. mult. & \#CG  & \#Exp.  & \#Prop.  \\ 
\hline
\multirow{3}{*}{fixed}
  & 2.0 & 615 & 1 & 306 & 1\\
  & 1.9 & 630 & 2 & 313 & 1\\
  & 1.95 & 640 & 2 & 318 & 1\\
\hline
 \multirow{3}{*}{grgr-optapprox}
  & 2.0 & 651 & 3 & 323 & 1\\
  & 1.9 & 663 & 1 & 330 & 1\\
  & 1.95 & 668 & 2 & 332 & 1\\
\hline
\multirow{3}{*}{gfgr-optapprox}
  & 2.0 & 621 & 3 & 308 & 1\\
  & 1.9 & 629 & 1 & 313 & 1\\
  & 1.95 & 638 & 2 & 317 & 1\\
\hline
 \multirow{3}{*}{grgr-opt}
  & 1.4 & 153 & 22 & 43 & 1\\
  & 1.6 & 169 & 26 & 47 & 1\\
  & 1.0 & 215 & 20 & 64 & 2\\
\hline
 \multirow{3}{*}{gfgr-opt}
  & 0.8 & 136 & 9 & 41 & 3\\
  & 1.4 & 173 & 25 & 48 & 3\\
  & 1.6 & 184 & 25 & 52 & 2\\
\hline
 \multirow{3}{*}{gfgf-opt}
  & 2.0 & 134 & 11 & 38 & 8\\
  & 1.8 & 156 & 18 & 44 & 5\\
  & 1.6 & 170 & 27 & 47 & 1\\
\hline
 projcg & - & 133 & 13 & 58 & 3 \\
\hline
\end{tabular}
\caption{Classification problem, Diabetes dataset: Comparison of expansion strategies by the number of SMALBE outer iterations, overall Hessian multiplications, CG, Expansion, and proportioning steps.}
\label{tab:diabetes}
\end{table}

\begin{table}
\centering
\begin{tabular}{| *{6}{c|}}
\hline
exp. type & $\alpha_u$ & \#Hess. mult. & \#CG  & \#Exp.  & \#Prop.  \\ 
\hline
\multirow{3}{*}{fixed}
  & 2 & 380 & 14 & 182 & 1\\
  & 1.9 & 381 & 21 & 179 & 1\\
  & 1.4 & 384 & 26 & 178 & 1\\
\hline
 \multirow{3}{*}{grgr-optapprox}
  & 1.95 & 393 & 9 & 191 & 1\\
  & 2 & 406 & 14 & 195 & 1\\
  & 1.9 & 420 & 14 & 202 & 1\\
\hline
\multirow{3}{*}{gfgr-optapprox}
  & 1.95 & 372 & 10 & 180 & 1\\
  & 1.4 & 397 & 33 & 181 & 1\\
  & 2 & 398 & 14 & 191 & 1\\
\hline
 \multirow{3}{*}{grgr-opt}
  & 0.8 & 141 & 8 & 43 & 3\\
  & 0.6 & 174 & 11 & 53 & 3\\
  & 2 & 180 & 23 & 50 & 6\\
\hline
 \multirow{3}{*}{gfgr-opt}
  & 1.6 & 160 & 30 & 42 & 3\\
  & 1.8 & 169 & 30 & 44 & 6\\
  & 1.2 & 180 & 18 & 52 & 5\\
\hline
 \multirow{3}{*}{gfgf-opt}
  & 2 & 113 & 22 & 29 & 3\\
  & 1.9 & 142 & 26 & 37 & 4\\
  & 1.95 & 173 & 28 & 46 & 6\\
\hline
 projcg & - & 125 & 14 & 54 & 2 \\
\hline
\end{tabular}
\caption{Classification problem, Ionosphere dataset: Comparison of expansion strategies by the number of SMALBE outer iterations, overall Hessian multiplications, CG, Expansion, and proportioning steps.}
\label{tab:ionosphere}
\end{table}

For each dataset, the default expansion with the \emph{fixed} step length achieves the best performance for $\alpha_u \approx 1.9$. The best \emph{fixed} strategy converged in $195$, $615$, and $380$  Hessian multiplications for Australian, Diabetes, and Ionosphere datasets, respectively.

As in the previous benchmark, the \emph{optapprox} strategy is quite close to the \emph{fixed} step length.

For the Australian dataset, the best \emph{opt} variant was \emph{gfgr} for $\alpha_u = 1.4$ which was closely followed by \emph{gfgf} for $\alpha_u = 2.0$.
Moreover, \emph{gfgf} with the same value of $\alpha_u$ was the best variant for the other datasets.
Comparing, for the best \emph{opt} strategy, the number of Hessian multiplication with \emph{fixed} we have $101$ (reduction of $48\%$) for Australian, $134$ (reduction of $78\%$) for Diabetes, $113$ (reduction of $70\%$) for Ionosphere.

The \emph{projcg} strategy performed the best for both Australian and Diabetes dataset.
It was the second-best from all results for the Ionosphere dataset.
Againg comparing the number of Hessian multiplication with \emph{fixed} we have $83$ (reduction of $57\%$) for Australian, $133$ (reduction of $78\%$) for Diabetes, $125$ (reduction of $67\%$) for Ionosphere.

\section{Conclusion}
\label{sec:conc}
The part of MPRGP algorithm dealing with the expansion of the active set was investigated.
The numerical results show that there can be a large number of expansions necessary to achieve convergence while using the default \emph{fixed} step length expansion.
Since the expansion is the most expensive part of the algorithm, we presented three alternative schemes for the expansion aimed at reducing the number of expansions.

The first two schemes are based on the optimal step length for the decrease of the quadratic cost function. 
Moreover, we can use any combination of the free gradient $\Vg^f$ and the reduced free gradient $\Vg^r$ for the computation of the step length and its direction.
The \emph{opt} step length is expensive due to additional Hessian multiplication.
Therefore, in \emph{optapprox} step length, we approximate the \emph{opt} step length by replacing multiplication by the Hessian with multiplication by the largest eigenvalue of the Hessian.

The third scheme replaces the expansion by a full CG step projected back into the feasible set.

We compared the expansion schemes on two benchmarks.
The first one was a 3D linear elasticity contact problem employing TFETI and SMALBE as the outer solver.
The second benchmark was classification problems on three datasets using SVMs.

The benchmarks showed that the optimal $\alpha_u$ for the \emph{fixed} strategy is about $1.9$.
Using \emph{optapprox} step length gives similar results to \emph{fixed}.
For most $\alpha_u$ the \emph{opt} strategies outperform \emph{fixed}.
The drawback of \emph{opt} is that there are relatively large jumps in the number of Hessian multiplications depending on $\alpha_u$.
It seems that it is best to use $\Vg^f$ for both step length computation and direction,
because using $\Vg^r$ necessitates the approximation of maximal eigenvalue of the Hessian.
Such approximation is typically made using the power iteration and usually needs $20-50$ extra Hessian multiplications.

The effectiveness of the alternative approaches is relatively low for the contact problem ($8\%$ reduction in the number of Hessian multiplications for both \emph{gfgf-opt} and \emph{projcg}).
However, in classification problems, due to a large number of the expansion steps needed by \emph{fixed}, the effectiveness is excellent.
The reduction of the number of Hessian multiplications ranges from $48 \%$ to $78 \%$ for \emph{opt} and from $57 \%$ to $78 \%$ for \emph{projcg} (excluding computation of the maximal eigenvalue).

The \emph{projcg} strategy consistently performed the best or was very close.
Moreover, it does not need an estimate of the maximal eigenvalue, nor a user-selected $\alpha_u$.
Additionally, the implementation of the algorithm is simplified.
Therefore, we recommend using the projected CG step in place of the standard \emph{fixed} step length expansion.

\section*{Acknowledgements}

The authors acknowledge the support of the Czech Science Foundation (GACR) project no. 17-22615S ...

\bibliographystyle{newsaxe}
\bibliography{main}

\appendix
\section{Full Numerical Results}
\label{sec:appendix}
\begin{table}
\centering
\begin{tabular}{|c|c|c|c|c|c|c|}
\hline
exp. type & $\alpha_u$ & outer it. & \#Hess. mult. & \#CG & \#Exp. & \#Prop. \\ 
\hline
\hline
fixed & 0.2 & 11 & 1026 & 167 & 423 & 2\\ \hline
fixed & 0.4 & 11 & 791 & 148 & 315 & 2\\ \hline
fixed & 0.6 & 11 & 667 & 154 & 250 & 2\\ \hline
fixed & 0.8 & 13 & 578 & 160 & 201 & 3\\ \hline
fixed & 1.0 & 8 & 463 & 117 & 168 & 2\\ \hline
fixed & 1.2 & 13 & 468 & 152 & 150 & 3\\ \hline
fixed & 1.4 & 11 & 427 & 150 & 132 & 2\\ \hline
fixed & 1.6 & 13 & 432 & 156 & 130 & 3\\ \hline
fixed & 1.8 & 10 & 318 & 158 & 74 & 2\\ \hline
fixed & 1.9 & 10 & 323 & 144 & 83 & 3\\ \hline
fixed & 1.95 & 7 & 333 & 70 & 127 & 2\\ \hline
fixed & 2.0 & 7 & 388 & 99 & 140 & 2\\ \hline\hline

grgr-optapprox & 0.2 & 8 & 1125 & 135 & 490 & 2\\ \hline
grgr-optapprox & 0.4 & 11 & 855 & 160 & 341 & 2\\ \hline
grgr-optapprox & 0.6 & 11 & 692 & 155 & 262 & 2\\ \hline
grgr-optapprox & 0.8 & 8 & 531 & 119 & 201 & 2\\ \hline
grgr-optapprox & 1.0 & 16 & 657 & 136 & 251 & 3\\ \hline
grgr-optapprox & 1.2 & 11 & 480 & 147 & 160 & 2\\ \hline
grgr-optapprox & 1.4 & 11 & 442 & 149 & 140 & 2\\ \hline
grgr-optapprox & 1.6 & 13 & 460 & 164 & 140 & 3\\ \hline
grgr-optapprox & 1.8 & 10 & 321 & 145 & 82 & 2\\ \hline
grgr-optapprox & 1.9 & 10 & 311 & 162 & 68 & 3\\ \hline
grgr-optapprox & 1.95 & 10 & 336 & 125 & 99 & 3\\ \hline
grgr-optapprox & 2.0 & 7 & 355 & 76 & 135 & 2\\ \hline\hline

gfgr-optapprox & 0.2 & 11 & 1026 & 167 & 423 & 2\\ \hline
gfgr-optapprox & 0.4 & 11 & 791 & 148 & 315 & 2\\ \hline
gfgr-optapprox & 0.6 & 11 & 667 & 154 & 250 & 2\\ \hline
gfgr-optapprox & 0.8 & 13 & 578 & 160 & 201 & 3\\ \hline
gfgr-optapprox & 1.0 & 8 & 464 & 118 & 168 & 2\\ \hline
gfgr-optapprox & 1.2 & 13 & 468 & 152 & 150 & 3\\ \hline
gfgr-optapprox & 1.4 & 11 & 427 & 150 & 132 & 2\\ \hline
gfgr-optapprox & 1.6 & 13 & 432 & 156 & 130 & 3\\ \hline
gfgr-optapprox & 1.8 & 10 & 317 & 159 & 73 & 2\\ \hline
gfgr-optapprox & 1.9 & 10 & 325 & 140 & 86 & 3\\ \hline
gfgr-optapprox & 1.95 & 8 & 333 & 95 & 114 & 2\\ \hline
gfgr-optapprox & 2.0 & 8 & 405 & 105 & 145 & 2\\ \hline
\end{tabular}
\caption{Linear elasticity contact problem: Comparison of expansion strategies by the number of SMALBE outer iterations, overall Hessian multiplications, CG, Expansion, and proportioning steps. Part 1/2.}
\end{table}

\begin{table}
\centering
\begin{tabular}{|c|c|c|c|c|c|c|}
\hline
exp. type & $\alpha_u$ & outer iters & \#Hess. mult. & \#CG  & \#Exp.  & \#Prop.  \\ 
\hline
\hline
grgr-opt & 0.2 & 8 & 315 & 113 & 64 & 2\\ \hline
grgr-opt & 0.4 & 8 & 359 & 124 & 75 & 2\\ \hline
grgr-opt & 0.6 & 9 & 316 & 149 & 52 & 2\\ \hline
grgr-opt & 0.8 & 9 & 297 & 132 & 51 & 3\\ \hline
grgr-opt & 1.0 & 7 & 401 & 94 & 99 & 3\\ \hline
grgr-opt & 1.2 & 11 & 342 & 186 & 47 & 4\\ \hline
grgr-opt & 1.4 & 9 & 405 & 154 & 80 & 2\\ \hline
grgr-opt & 1.6 & 8 & 292 & 120 & 54 & 2\\ \hline
grgr-opt & 1.8 & 17 & 352 & 143 & 63 & 3\\ \hline
grgr-opt & 1.9 & 17 & 717 & 115 & 194 & 3\\ \hline
grgr-opt & 1.95 & 21 & 408 & 174 & 70 & 3\\ \hline
grgr-opt & 2.0 & 13 & 382 & 150 & 72 & 3\\ \hline\hline

gfgr-opt & 0.2 & 8 & 323 & 142 & 57 & 2\\ \hline
gfgr-opt & 0.4 & 8 & 316 & 111 & 65 & 2\\ \hline
gfgr-opt & 0.6 & 12 & 286 & 139 & 44 & 3\\ \hline
gfgr-opt & 0.8 & 9 & 327 & 123 & 64 & 3\\ \hline
gfgr-opt & 1.0 & 10 & 702 & 161 & 176 & 3\\ \hline
gfgr-opt & 1.2 & 10 & 326 & 154 & 53 & 3\\ \hline
gfgr-opt & 1.4 & 9 & 289 & 154 & 41 & 3\\ \hline
gfgr-opt & 1.6 & 9 & 269 & 129 & 43 & 2\\ \hline
gfgr-opt & 1.8 & 15 & 379 & 185 & 58 & 5\\ \hline
gfgr-opt & 1.9 & 13 & 618 & 179 & 141 & 3\\ \hline
gfgr-opt & 1.95 & 13 & 290 & 127 & 49 & 3\\ \hline
gfgr-opt & 2.0 & 14 & 484 & 181 & 95 & 4\\ \hline\hline

gfgf-opt & 0.2 & 9 & 311 & 135 & 55 & 2\\ \hline
gfgf-opt & 0.4 & 10 & 353 & 122 & 73 & 2\\ \hline
gfgf-opt & 0.6 & 13 & 326 & 171 & 46 & 4\\ \hline
gfgf-opt & 0.8 & 9 & 394 & 109 & 91 & 3\\ \hline
gfgf-opt & 1.0 & 11 & 297 & 171 & 37 & 4\\ \hline
gfgf-opt & 1.2 & 10 & 292 & 151 & 43 & 2\\ \hline
gfgf-opt & 1.4 & 9 & 299 & 152 & 45 & 3\\ \hline
gfgf-opt & 1.6 & 15 & 323 & 187 & 39 & 4\\ \hline
gfgf-opt & 1.8 & 17 & 454 & 134 & 100 & 3\\ \hline
gfgf-opt & 1.9 & 18 & 334 & 201 & 37 & 4\\ \hline
gfgf-opt & 1.95 & 18 & 337 & 192 & 41 & 4\\ \hline
gfgf-opt & 2.0 & 16 & 384 & 208 & 52 & 4\\ \hline\hline

projcg & - & 10 & 292 & 171 & 53 & 5\\ \hline
\end{tabular}
\caption{Linear elasticity contact problem: Comparison of expansion strategies by the number of SMALBE outer iterations, overall Hessian multiplications, CG, Expansion, and proportioning steps. Part 2/2.}
\end{table}

\begin{table}
\centering
\begin{tabular}{|c|c|c|c|c|c|}
\hline
exp. type & $\alpha_u$ & \#Hess. mult. & \#CG  & \#Exp.  & \#Prop.  \\ 
\hline
\hline
grgr-fixed & 0.2 & 441 & 31 & 204 & 1\\ \hline
grgr-fixed & 0.4 & 355 & 29 & 162 & 1\\ \hline
grgr-fixed & 0.6 & 309 & 29 & 139 & 1\\ \hline
grgr-fixed & 0.8 & 298 & 30 & 133 & 1\\ \hline
grgr-fixed & 1.0 & 261 & 21 & 119 & 1\\ \hline
grgr-fixed & 1.2 & 248 & 22 & 112 & 1\\ \hline
grgr-fixed & 1.4 & 232 & 22 & 104 & 1\\ \hline
grgr-fixed & 1.6 & 227 & 24 & 100 & 2\\ \hline
grgr-fixed & 1.8 & 205 & 9 & 97 & 1\\ \hline
grgr-fixed & 1.9 & 195 & 8 & 92 & 2\\ \hline
grgr-fixed & 1.95 & 195 & 8 & 92 & 2\\ \hline
grgr-fixed & 2.0 & 202 & 10 & 95 & 1\\ \hline
\hline
grgr-optapprox & 0.2 & 438 & 22 & 207 & 1\\ \hline
grgr-optapprox & 0.4 & 362 & 28 & 166 & 1\\ \hline
grgr-optapprox & 0.6 & 319 & 25 & 146 & 1\\ \hline
grgr-optapprox & 0.8 & 306 & 34 & 135 & 1\\ \hline
grgr-optapprox & 1.0 & 281 & 23 & 128 & 1\\ \hline
grgr-optapprox & 1.2 & 259 & 21 & 118 & 1\\ \hline
grgr-optapprox & 1.4 & 244 & 20 & 111 & 1\\ \hline
grgr-optapprox & 1.6 & 239 & 21 & 108 & 1\\ \hline
grgr-optapprox & 1.8 & 220 & 14 & 102 & 1\\ \hline
grgr-optapprox & 1.9 & 195 & 9 & 92 & 1\\ \hline
grgr-optapprox & 1.95 & 211 & 7 & 101 & 1\\ \hline
grgr-optapprox & 2.0 & 224 & 1 & 110 & 2\\ \hline
\hline
gfgr-optapprox & 0.2 & 437 & 33 & 201 & 1\\ \hline
gfgr-optapprox & 0.4 & 344 & 26 & 158 & 1\\ \hline
gfgr-optapprox & 0.6 & 306 & 27 & 138 & 2\\ \hline
gfgr-optapprox & 0.8 & 293 & 27 & 132 & 1\\ \hline
gfgr-optapprox & 1.0 & 251 & 19 & 115 & 1\\ \hline
gfgr-optapprox & 1.2 & 247 & 21 & 112 & 1\\ \hline
gfgr-optapprox & 1.4 & 227 & 23 & 101 & 1\\ \hline
gfgr-optapprox & 1.6 & 225 & 29 & 97 & 1\\ \hline
gfgr-optapprox & 1.8 & 203 & 13 & 94 & 1\\ \hline
gfgr-optapprox & 1.9 & 198 & 16 & 90 & 1\\ \hline
gfgr-optapprox & 1.95 & 187 & 5 & 90 & 1\\ \hline
gfgr-optapprox & 2.0 & 199 & 5 & 96 & 1\\ \hline
\end{tabular}
\caption{Classification problem, Australian dataset: Comparison of expansion strategies by the number of overall Hessian multiplications, CG, Expansion, and proportioning steps. Part 1/2.}
\end{table}
\begin{table}
\centering
\begin{tabular}{|c|c|c|c|c|c|}
\hline
exp. type & $\alpha_u$ & \#Hess. mult. & \#CG  & \#Exp.  & \#Prop.  \\ 
\hline
\hline
grgr-opt & 0.2 & 263 & 24 & 79 & 1\\ \hline
grgr-opt & 0.4 & 202 & 20 & 60 & 1\\ \hline
grgr-opt & 0.6 & 179 & 18 & 53 & 1\\ \hline
grgr-opt & 0.8 & 163 & 16 & 48 & 2\\ \hline
grgr-opt & 1.0 & 188 & 14 & 57 & 2\\ \hline
grgr-opt & 1.2 & 187 & 21 & 54 & 3\\ \hline
grgr-opt & 1.4 & 153 & 27 & 41 & 2\\ \hline
grgr-opt & 1.6 & 177 & 24 & 50 & 2\\ \hline
grgr-opt & 1.8 & 166 & 27 & 45 & 3\\ \hline
grgr-opt & 1.9 & 148 & 28 & 39 & 2\\ \hline
grgr-opt & 1.95 & 300 & 48 & 82 & 5\\ \hline
grgr-opt & 2.0 & 269 & 45 & 73 & 4\\ \hline
\hline
gfgr-opt & 0.2 & 274 & 26 & 82 & 1\\ \hline
gfgr-opt & 0.4 & 177 & 16 & 53 & 1\\ \hline
gfgr-opt & 0.6 & 176 & 21 & 51 & 1\\ \hline
gfgr-opt & 0.8 & 161 & 18 & 47 & 1\\ \hline
gfgr-opt & 1.0 & 135 & 15 & 39 & 2\\ \hline
gfgr-opt & 1.2 & 165 & 27 & 45 & 2\\ \hline
gfgr-opt & 1.4 & 101 & 20 & 26 & 2\\ \hline
gfgr-opt & 1.6 & 194 & 28 & 54 & 3\\ \hline
gfgr-opt & 1.8 & 140 & 25 & 37 & 3\\ \hline
gfgr-opt & 1.9 & 173 & 33 & 45 & 4\\ \hline
gfgr-opt & 1.95 & 244 & 41 & 65 & 7\\ \hline
gfgr-opt & 2.0 & 187 & 35 & 49 & 4\\ \hline
\hline
gfgf-opt & 0.2 & 244 & 20 & 74 & 1\\ \hline
gfgf-opt & 0.4 & 192 & 24 & 55 & 2\\ \hline
gfgf-opt & 0.6 & 141 & 16 & 41 & 1\\ \hline
gfgf-opt & 0.8 & 145 & 17 & 42 & 1\\ \hline
gfgf-opt & 1.0 & 156 & 22 & 44 & 1\\ \hline
gfgf-opt & 1.2 & 159 & 28 & 43 & 1\\ \hline
gfgf-opt & 1.4 & 129 & 28 & 33 & 1\\ \hline
gfgf-opt & 1.6 & 207 & 39 & 55 & 2\\ \hline
gfgf-opt & 1.8 & 165 & 36 & 42 & 2\\ \hline
gfgf-opt & 1.9 & 124 & 28 & 31 & 2\\ \hline
gfgf-opt & 1.95 & 151 & 34 & 38 & 2\\ \hline
gfgf-opt & 2.0 & 107 & 23 & 27 & 2\\ \hline
\hline
projcg & - & 83 & 16 & 32 & 2\\ \hline
\end{tabular}
\caption{Classification problem, Australian dataset: Comparison of expansion strategies by the number of overall Hessian multiplications, CG, Expansion, and proportioning steps. Part 2/2.}
\end{table}

\begin{table}
\centering
\begin{tabular}{|c|c|c|c|c|c|}
\hline
exp. type & $\alpha_u$ & \#Hess. mult. & \#CG  & \#Exp.  & \#Prop.  \\ 
\hline
\hline
grgr-fixed & 0.2 & 1065 & 65 & 499 & 1\\ \hline
grgr-fixed & 0.4 & 976 & 52 & 461 & 1\\ \hline
grgr-fixed & 0.6 & 887 & 51 & 417 & 1\\ \hline
grgr-fixed & 0.8 & 851 & 31 & 409 & 1\\ \hline
grgr-fixed & 1.0 & 800 & 14 & 392 & 1\\ \hline
grgr-fixed & 1.2 & 765 & 9 & 377 & 1\\ \hline
grgr-fixed & 1.4 & 697 & 5 & 345 & 1\\ \hline
grgr-fixed & 1.6 & 653 & 3 & 324 & 1\\ \hline
grgr-fixed & 1.8 & 648 & 4 & 321 & 1\\ \hline
grgr-fixed & 1.9 & 630 & 2 & 313 & 1\\ \hline
grgr-fixed & 1.95 & 640 & 2 & 318 & 1\\ \hline
grgr-fixed & 2.0 & 615 & 1 & 306 & 1\\ \hline
\hline
grgr-optapprox & 0.2 & 1080 & 66 & 506 & 1\\ \hline
grgr-optapprox & 0.4 & 988 & 52 & 467 & 1\\ \hline
grgr-optapprox & 0.6 & 909 & 47 & 430 & 1\\ \hline
grgr-optapprox & 0.8 & 890 & 32 & 428 & 1\\ \hline
grgr-optapprox & 1.0 & 823 & 15 & 403 & 1\\ \hline
grgr-optapprox & 1.2 & 804 & 14 & 394 & 1\\ \hline
grgr-optapprox & 1.4 & 737 & 5 & 365 & 1\\ \hline
grgr-optapprox & 1.6 & 687 & 3 & 341 & 1\\ \hline
grgr-optapprox & 1.8 & 684 & 2 & 340 & 1\\ \hline
grgr-optapprox & 1.9 & 663 & 1 & 330 & 1\\ \hline
grgr-optapprox & 1.95 & 668 & 2 & 332 & 1\\ \hline
grgr-optapprox & 2.0 & 651 & 3 & 323 & 1\\ \hline
\hline
gfgr-optapprox & 0.2 & 1070 & 66 & 501 & 1\\ \hline
gfgr-optapprox & 0.4 & 976 & 52 & 461 & 1\\ \hline
gfgr-optapprox & 0.6 & 888 & 50 & 418 & 1\\ \hline
gfgr-optapprox & 0.8 & 851 & 31 & 409 & 1\\ \hline
gfgr-optapprox & 1.0 & 802 & 16 & 392 & 1\\ \hline
gfgr-optapprox & 1.2 & 765 & 9 & 377 & 1\\ \hline
gfgr-optapprox & 1.4 & 697 & 5 & 345 & 1\\ \hline
gfgr-optapprox & 1.6 & 653 & 3 & 324 & 1\\ \hline
gfgr-optapprox & 1.8 & 648 & 4 & 321 & 1\\ \hline
gfgr-optapprox & 1.9 & 629 & 1 & 313 & 1\\ \hline
gfgr-optapprox & 1.95 & 638 & 2 & 317 & 1\\ \hline
gfgr-optapprox & 2.0 & 621 & 3 & 308 & 1\\ \hline
\end{tabular}
\caption{Classification problem, Diabetes dataset: Comparison of expansion strategies by the number of overall Hessian multiplications, CG, Expansion, and proportioning steps. Part 1/2.}
\end{table}
\begin{table}
\centering
\begin{tabular}{|c|c|c|c|c|c|}
\hline
exp. type & $\alpha_u$ & \#Hess. mult. & \#CG  & \#Exp.  & \#Prop.  \\ 
\hline
\hline
grgr-opt & 0.2 & 311 & 13 & 98 & 3\\ \hline
grgr-opt & 0.4 & 318 & 20 & 98 & 3\\ \hline
grgr-opt & 0.6 & 284 & 10 & 89 & 6\\ \hline
grgr-opt & 0.8 & 237 & 20 & 71 & 3\\ \hline
grgr-opt & 1.0 & 215 & 20 & 64 & 2\\ \hline
grgr-opt & 1.2 & 254 & 29 & 74 & 2\\ \hline
grgr-opt & 1.4 & 153 & 22 & 43 & 1\\ \hline
grgr-opt & 1.6 & 169 & 26 & 47 & 1\\ \hline
grgr-opt & 1.8 & 240 & 30 & 68 & 5\\ \hline
grgr-opt & 1.9 & 290 & 34 & 82 & 9\\ \hline
grgr-opt & 1.95 & 273 & 23 & 80 & 9\\ \hline
grgr-opt & 2.0 & 309 & 31 & 89 & 10\\ \hline
\hline
gfgr-opt & 0.2 & 256 & 13 & 80 & 2\\ \hline
gfgr-opt & 0.4 & 230 & 12 & 71 & 4\\ \hline
gfgr-opt & 0.6 & 311 & 36 & 91 & 1\\ \hline
gfgr-opt & 0.8 & 136 & 9 & 41 & 3\\ \hline
gfgr-opt & 1.0 & 308 & 38 & 89 & 2\\ \hline
gfgr-opt & 1.2 & 244 & 43 & 66 & 2\\ \hline
gfgr-opt & 1.4 & 173 & 25 & 48 & 3\\ \hline
gfgr-opt & 1.6 & 184 & 25 & 52 & 2\\ \hline
gfgr-opt & 1.8 & 211 & 34 & 57 & 5\\ \hline
gfgr-opt & 1.9 & 196 & 28 & 54 & 5\\ \hline
gfgr-opt & 1.95 & 320 & 28 & 91 & 18\\ \hline
gfgr-opt & 2.0 & 280 & 34 & 77 & 14\\ \hline
\hline
gfgf-opt & 0.2 & 289 & 16 & 90 & 2\\ \hline
gfgf-opt & 0.4 & 233 & 10 & 73 & 3\\ \hline
gfgf-opt & 0.6 & 361 & 20 & 111 & 7\\ \hline
gfgf-opt & 0.8 & 308 & 31 & 90 & 6\\ \hline
gfgf-opt & 1.0 & 299 & 34 & 87 & 3\\ \hline
gfgf-opt & 1.2 & 194 & 24 & 55 & 4\\ \hline
gfgf-opt & 1.4 & 334 & 46 & 94 & 5\\ \hline
gfgf-opt & 1.6 & 170 & 27 & 47 & 1\\ \hline
gfgf-opt & 1.8 & 156 & 18 & 44 & 5\\ \hline
gfgf-opt & 1.9 & 253 & 32 & 71 & 7\\ \hline
gfgf-opt & 1.95 & 250 & 18 & 72 & 15\\ \hline
gfgf-opt & 2.0 & 134 & 11 & 38 & 8\\ \hline
\hline
projcg & - & 133 & 13 & 58 & 3\\ \hline
\end{tabular}
\caption{Classification problem, Diabetes dataset: Comparison of expansion strategies by the number of overall Hessian multiplications, CG, Expansion, and proportioning steps. Part 2/2.}
\end{table}

\begin{table}
\centering
\begin{tabular}{|c|c|c|c|c|c|}
\hline
exp. type & $\alpha_u$ & \#Hess. mult. & \#CG  & \#Exp.  & \#Prop.  \\ 
\hline
\hline
grgr-fixed & 0.2 & 565 & 53 & 255 & 1\\ \hline
grgr-fixed & 0.4 & 514 & 50 & 231 & 1\\ \hline
grgr-fixed & 0.6 & 428 & 26 & 200 & 1\\ \hline
grgr-fixed & 0.8 & 443 & 38 & 201 & 2\\ \hline
grgr-fixed & 1.0 & 444 & 48 & 197 & 1\\ \hline
grgr-fixed & 1.2 & 428 & 44 & 191 & 1\\ \hline
grgr-fixed & 1.4 & 384 & 26 & 178 & 1\\ \hline
grgr-fixed & 1.6 & 403 & 26 & 187 & 2\\ \hline
grgr-fixed & 1.8 & 437 & 37 & 199 & 1\\ \hline
grgr-fixed & 1.9 & 381 & 21 & 179 & 1\\ \hline
grgr-fixed & 1.95 & 393 & 15 & 188 & 1\\ \hline
grgr-fixed & 2.0 & 380 & 14 & 182 & 1\\ \hline
\hline
grgr-optapprox & 0.2 & 565 & 51 & 256 & 1\\ \hline
grgr-optapprox & 0.4 & 525 & 53 & 235 & 1\\ \hline
grgr-optapprox & 0.6 & 447 & 25 & 210 & 1\\ \hline
grgr-optapprox & 0.8 & 468 & 39 & 213 & 2\\ \hline
grgr-optapprox & 1.0 & 450 & 42 & 203 & 1\\ \hline
grgr-optapprox & 1.2 & 452 & 48 & 201 & 1\\ \hline
grgr-optapprox & 1.4 & 424 & 34 & 194 & 1\\ \hline
grgr-optapprox & 1.6 & 433 & 30 & 200 & 2\\ \hline
grgr-optapprox & 1.8 & 425 & 23 & 200 & 1\\ \hline
grgr-optapprox & 1.9 & 420 & 14 & 202 & 1\\ \hline
grgr-optapprox & 1.95 & 393 & 9 & 191 & 1\\ \hline
grgr-optapprox & 2.0 & 406 & 14 & 195 & 1\\ \hline
\hline
gfgr-optapprox & 0.2 & 560 & 52 & 253 & 1\\ \hline
gfgr-optapprox & 0.4 & 520 & 54 & 232 & 1\\ \hline
gfgr-optapprox & 0.6 & 428 & 22 & 202 & 1\\ \hline
gfgr-optapprox & 0.8 & 453 & 42 & 204 & 2\\ \hline
gfgr-optapprox & 1.0 & 432 & 42 & 194 & 1\\ \hline
gfgr-optapprox & 1.2 & 428 & 44 & 191 & 1\\ \hline
gfgr-optapprox & 1.4 & 397 & 33 & 181 & 1\\ \hline
gfgr-optapprox & 1.6 & 406 & 27 & 188 & 2\\ \hline
gfgr-optapprox & 1.8 & 430 & 30 & 199 & 1\\ \hline
gfgr-optapprox & 1.9 & 425 & 17 & 203 & 1\\ \hline
gfgr-optapprox & 1.95 & 372 & 10 & 180 & 1\\ \hline
gfgr-optapprox & 2.0 & 398 & 14 & 191 & 1\\ \hline
\end{tabular}
\caption{Classification problem, Ionosphere dataset: Comparison of expansion strategies by the number of overall Hessian multiplications, CG, Expansion, and proportioning steps. Part 1/2.}
\end{table}
\begin{table}
\centering
\begin{tabular}{|c|c|c|c|c|c|}
\hline
exp. type & $\alpha_u$ & \#Hess. mult. & \#CG  & \#Exp.  & \#Prop.  \\ 
\hline
\hline
grgr-opt & 0.2 & 264 & 25 & 79 & 1\\ \hline
grgr-opt & 0.4 & 264 & 17 & 81 & 3\\ \hline
grgr-opt & 0.6 & 174 & 11 & 53 & 3\\ \hline
grgr-opt & 0.8 & 141 & 8 & 43 & 3\\ \hline
grgr-opt & 1.0 & 193 & 15 & 58 & 3\\ \hline
grgr-opt & 1.2 & 226 & 21 & 67 & 3\\ \hline
grgr-opt & 1.4 & 280 & 39 & 78 & 6\\ \hline
grgr-opt & 1.6 & 187 & 32 & 50 & 4\\ \hline
grgr-opt & 1.8 & 187 & 26 & 52 & 4\\ \hline
grgr-opt & 1.9 & 233 & 32 & 65 & 5\\ \hline
grgr-opt & 1.95 & 339 & 27 & 99 & 14\\ \hline
grgr-opt & 2.0 & 180 & 23 & 50 & 6\\ \hline
\hline
gfgr-opt & 0.2 & 309 & 31 & 92 & 1\\ \hline
gfgr-opt & 0.4 & 211 & 22 & 62 & 2\\ \hline
gfgr-opt & 0.6 & 193 & 26 & 55 & 1\\ \hline
gfgr-opt & 0.8 & 198 & 31 & 55 & 1\\ \hline
gfgr-opt & 1.0 & 181 & 28 & 50 & 2\\ \hline
gfgr-opt & 1.2 & 180 & 18 & 52 & 5\\ \hline
gfgr-opt & 1.4 & 221 & 38 & 59 & 5\\ \hline
gfgr-opt & 1.6 & 160 & 30 & 42 & 3\\ \hline
gfgr-opt & 1.8 & 169 & 30 & 44 & 6\\ \hline
gfgr-opt & 1.9 & 2322 & 57 & 659 & 287\\ \hline
gfgr-opt & 1.95 & 233 & 25 & 64 & 15\\ \hline
gfgr-opt & 2.0 & 219 & 31 & 59 & 10\\ \hline
\hline
gfgf-opt & 0.2 & 256 & 26 & 76 & 1\\ \hline
gfgf-opt & 0.4 & 244 & 28 & 71 & 2\\ \hline
gfgf-opt & 0.6 & 206 & 16 & 62 & 3\\ \hline
gfgf-opt & 0.8 & 229 & 19 & 68 & 5\\ \hline
gfgf-opt & 1.0 & 200 & 30 & 56 & 1\\ \hline
gfgf-opt & 1.2 & 186 & 19 & 54 & 4\\ \hline
gfgf-opt & 1.4 & 200 & 30 & 55 & 4\\ \hline
gfgf-opt & 1.6 & 175 & 32 & 46 & 4\\ \hline
gfgf-opt & 1.8 & 216 & 29 & 58 & 12\\ \hline
gfgf-opt & 1.9 & 142 & 26 & 37 & 4\\ \hline
gfgf-opt & 1.95 & 173 & 28 & 46 & 6\\ \hline
gfgf-opt & 2.0 & 113 & 22 & 29 & 3\\ \hline
\hline
projcg & - & 125 & 14 & 54 & 2\\ \hline
\end{tabular}
\caption{Classification problem, Ionosphere dataset: Comparison of expansion strategies by the number of overall Hessian multiplications, CG, Expansion, and proportioning steps. Part 2/2.}
\end{table}
\end{document}